\documentclass[11pt]{amsart}

\usepackage{amsmath,amssymb,amsthm,amscd,amsfonts}
\usepackage{verbatim,psfrag}
\usepackage[all]{xy}
\usepackage{pinlabel}

\newcommand{\nc}{\newcommand}
\nc{\nt}{\newtheorem}
\nc{\bs}{\bigskip}
\nc{\dmo}{\DeclareMathOperator}

\nt*{main}{Main Theorem}
\nt{theorem}{Theorem}[section]
\nt{proposition}[theorem]{Proposition}
\nt{lemma}[theorem]{Lemma}
\nt{fact}[theorem]{Fact}
\nt{corollary}[theorem]{Corollary}
\nt{conjecture}[theorem]{Conjecture}
\nt{question}[theorem]{Question}
\nt{problem}[theorem]{Problem}

\dmo{\Aut}{Aut}
\dmo{\Mod}{Mod}
\dmo{\PMod}{PMod}
\dmo{\SMod}{SMod}
\dmo{\LMod}{LMod}
\dmo{\SHomeo}{SHomeo}
\dmo{\Homeo}{Homeo}
\dmo{\PHomeo}{PHomeo}
\dmo{\Teich}{Teich}
\dmo{\Sp}{Sp}
\dmo{\SL}{SL}
\dmo{\GL}{GL}

\nc{\Z}{\mathbb{Z}}
\nc{\R}{\mathbb{R}}
\nc{\N}{\mathbb{N}}
\nc{\C}{\mathbb{C}}

\dmo{\B}{B}
\dmo{\Br}{K}
\dmo{\PB}{PB}
\dmo{\even}{even}

\nc{\p}[1]{\medskip\paragraph{{\em #1}}}
\nc{\margin}[1]{\marginpar{\scriptsize #1}}
\nc{\bl}{ \begin{list}{$\cdot$}{
\setlength{\leftmargin}{.5in}
\setlength{\rightmargin}{.5in}
\setlength{\parsep}{0.5ex plus .2ex minus 0ex}
\setlength{\itemsep}{0.2ex plus 0.2ex minus 0ex}
}
}

\nc{\el}{\end{list}}

\title{The Birman--Hilden theory}

\begin{document}
	
\input{epsf.sty}

\author{Dan Margalit and Rebecca R. Winarski}

\address{Dan Margalit \\ School of Mathematics\\ Georgia Institute of Technology \\ 686 Cherry St. \\ Atlanta, GA 30332 \\  margalit@math.gatech.edu}

\thanks{The first author is supported by the National Science Foundation under Grant No. DMS - 1057874.}

\address{Rebecca R. Winarski\\ Department of Mathematical Sciences\\University of Wisconsin-Milwaukee\\ rebecca.winarski@gmail.com}
%\keywords{braid group, congruence subgroup}

%\subjclass[2000]{Primary: 20F36; Secondary: 57M07}

\begin{abstract}
In the 1970s Joan Birman and Hugh Hilden wrote several papers on the problem of relating the mapping class group of a surface to that of a cover.  We survey their work,  give an overview of the subsequent developments, and discuss open questions and new directions.
\end{abstract}

\maketitle

\vspace*{-2ex}

\section{Introduction}

In the early 1970s Joan Birman and Hugh Hilden wrote a series of now-classic papers on the interplay between mapping class groups and covering spaces.  The initial goal was to determine a presentation for the mapping class group of $S_2$, the closed surface of genus two (it was not until the late 1970s that Hatcher and Thurston \cite{HT} developed an approach for finding explicit presentations for mapping class groups).

The key innovation by Birman and Hilden is to relate the mapping class group $\Mod(S_2)$ to the mapping class group of $S_{0,6}$, a sphere with six marked points.  Presentations for $\Mod(S_{0,6})$  were already known since that group is closely related to a braid group.

%\begin{figure}[h!]
%\end{figure}

The two surfaces $S_2$ and $S_{0,6}$ are related by a two-fold branched covering map $S_2 \to S_{0,6}$:
\begin{center}
\includegraphics[scale=.65]{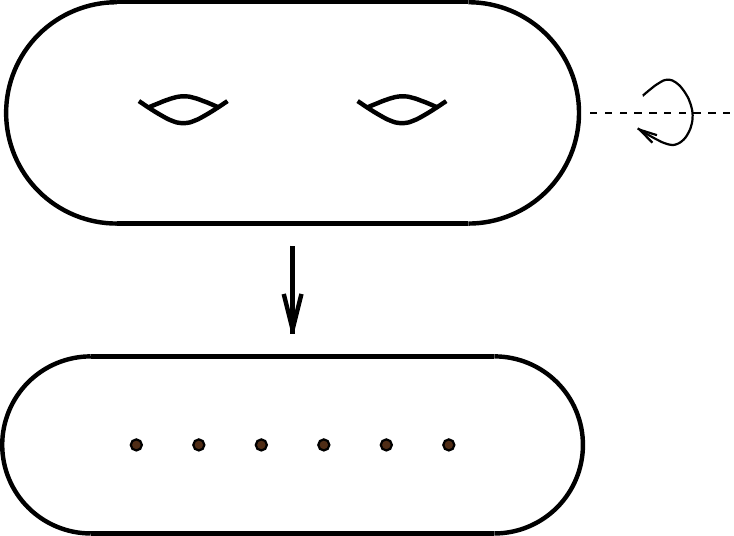}
\end{center}
The six marked points in the base are branch points.  The deck transformation is called the \emph{hyperelliptic involution} of $S_2$, and we denote it by $\iota$.  Every element of $\Mod(S_2)$ has a representative that commutes with $\iota$, and so it follows that there is a map 
\[
\Theta : \Mod(S_2) \to \Mod(S_{0,6}).
 \]
The kernel of $\Theta$ is the cyclic group of order two generated by (the homotopy class of) the involution $\iota$.  One can check that each generator for $\Mod(S_{0,6})$ lifts to $\Mod(S_2)$ and so $\Theta$ is surjective.  From this we have a short exact sequence
\[ 1 \to \langle \iota \rangle \to \Mod(S_2) \stackrel{\Theta}{\to} \Mod(S_{0,6}) \to 1,\]
and hence a presentation for $\Mod(S_{0,6})$ can be lifted to a presentation for $\Mod(S_2)$.

But wait---the map $\Theta$ is not a priori well defined!  The problem is that elements of $\Mod(S_2)$ are only defined up to isotopy, and these isotopies are not required to respect the hyperelliptic involution.  The first paper by Birman and Hilden proves that in fact all isotopies can be chosen to respect the involution.  Birman and Hilden quickly realized that the theory initiated in that first paper can be generalized in various ways, and they wrote a series of paper on the subject, culminating in the paper \emph{On Isotopies of Homeomorphisms of Riemann Surfaces} \cite{BH73}, published in \emph{Annals of Mathematics} in 1973.

In the remainder of this article, we will discuss the history of the Birman--Hilden theory, including generalizations by MacLachlan--Harvey and the second author of this article, we will give several applications, explain three proofs, and discuss various open questions and new directions in the theory.  As we will see, the Birman--Hilden theory has had influence on many areas of mathematics, from low-dimensional topology, to geometric group theory, to representation theory, to algebraic geometry and more, and it continues to produce interesting open problems and research directions.

\bigskip

\noindent\emph{The other article by Birman and Hilden.} Before getting on with our main business, we would be remiss not to mention the other paper by Birman and Hilden \cite{BH75}, the 1975 paper \emph{Heegaard splittings of branched coverings of $S^3$}, published in \emph{Transactions of the American Mathematical Society} (there is also the corresponding research announcement, \emph{The homeomorphism problem for $S^3$}, published two years earlier \cite{BH3mfd73}).  In this paper, Birman and Hilden discuss the relationship between branched covers and Heegaard splittings of 3-manifolds.  Their results cover a lot of territory.  For instance: 
\begin{itemize}
\item they prove that every closed, orientable 3-manifold of Heegaard genus 2 is a two-fold branched covering space of $S^3$ branched over a 3-bridge knot or link;
\item they give an algorithm for determining if a Heegaard splitting of genus two represents $S^3$;
\item they prove that any simply connected two-fold cover of $S^3$ branched over the closure of a braid on three strands is itself $S^3$; and
\item they disprove a conjecture of Haken that among all simply connected 3-manifolds, and among all group presentations for their fundamental groups arising from their Heegaard splittings, the presentations for $\pi_1(S^3)$ have a certain nice property.
\end{itemize}
While this paper has also been influential and well-cited, and in fact relies on their work on surfaces, we will restrict our focus in this article to the work of Birman and Hilden on mapping class groups.

\section{Statements of the main theorem}\label{state}

Let $p: S \to X$ be a covering map of surfaces, possibly branched, possibly with boundary.  We say that $f : S \to S$ is \emph{fiber preserving} if for each $x \in X$ there is a $y \in X$ so that $f(p^{-1}(x)) = p^{-1}(y)$; in other words, as the terminology suggests, $f$ takes fibers to fibers.  

Given two homotopic fiber-preserving homeomorphisms of $S$, we can ask if they are homotopic through fiber-preserving homeomorphisms.  If the answer is yes for all such pairs of homeomorphisms, we say that the covering map $p$ has the \emph{Birman--Hilden property}.  An equivalent formulation of the Birman--Hilden property is: whenever a fiber-preserving homeomorphism is homotopic to the identity, it is homotopic to the identity through fiber-preserving homeomorphisms.

We are now ready to state the main theorems of the Birman--Hilden theory.  There are several versions, proved over the years by various authors, each generalizing the previous.  The first version is the one that appears in the aforementioned 1973 \emph{Annals of Mathematics} paper by Birman and Hilden and also in the accompanying research announcement \emph{Isotopies of Homeomorphisms of Riemann surfaces} \cite{BHannouncement}.   Throughout, we will say that a surface is \emph{hyperbolic} if its Euler characteristic is negative.

\begin{theorem}[Birman--Hilden]\label{bh}
Let $p : S \to X$ be a finite-sheeted regular branched covering map where $S$ is a hyperbolic surface.  Assume that $p$ is either unbranched or is solvable.  Then $p$ has the Birman--Hilden property.
\end{theorem}

If we apply Theorem~\ref{bh} to the branched covering map $S_2 \to S_{0,6}$ described earlier, then it exactly says that the map $\Theta : \Mod(S_2) \to \Mod(S_{0,6})$ is well defined.

It is worthwhile to compare our Theorem~\ref{bh} to what is actually stated by Birman and Hilden.  In their paper, they state two theorems, each of which is a special case of Theorem~\ref{bh}.  Their Theorem~1 treats the case of regular covers where each deck transformation fixes each preimage of each branch point in $X$.  This clearly takes care of the case of unbranched covers, and also the case of certain solvable branched covers (on one hand  a finite group of homeomorphisms of a surface that fixes a point must be a subgroup of a dihedral group, and on the other hand there are solvable---even cyclic---branched covers that do not satisfy the condition of Theorem~1).  Birman and Hilden's Theorem~2 deals with the general case of solvable covers, which includes some unbranched covers.

In early 1973 MacLachlan and Harvey \cite{MH} published a paper called \emph{On Mapping Class Groups and Covering Spaces}, in which they give the following generalization of Theorem~\ref{bh}.

\begin{theorem}[MacLachlan--Harvey]\label{mh}
Let $p : S \to X$ be a finite-sheeted regular branched covering map where $S$ is a hyperbolic surface.  Then $p$ has the Birman--Hilden property.
\end{theorem}

MacLachlan and Harvey's work was contemporaneous with the work of Birman and Hilden cited in Theorem~\ref{bh}, and was subsequent to the original paper by Birman and Hilden on the hyperelliptic case.  Their approach is completely different, and is framed in terms of Teichm\"uller theory.  

The 2014 Ph.D. thesis of the second author of this article is a further generalization \cite{winarski}.  For the statement, a preimage of a branch point is \emph{unramified} if some small neighborhood is mapped injectively under the covering map, and a cover is \emph{fully ramified} if no branch point has an unramified preimage.

\begin{theorem}[Winarski]\label{win}
Let $p: S \to X$ be a finite-sheeted branched covering map where $S$ is a hyperbolic surface, and suppose that $p$ is fully ramified.  Then $p$ has the Birman--Hilden property.  
\end{theorem}

Note that all regular covers are fully ramified and also that all unbranched covers are fully ramified.  Thus Theorem~\ref{win} indeed implies Theorems~\ref{bh} and~\ref{mh}.  In Section 2.3 of her paper, Winarski gives a general construction of irregular branched covers that are fully ramified.  Thus there are many examples of covering spaces that satisfy the hypotheses of Theorem~\ref{win} but not those of Theorem~\ref{mh}.

We will briefly remark on the assumption that $S$ is hyperbolic.  It is not hard to construct counterexamples in the other cases.  For instance suppose $S$ is the torus $T^2$ and $p : S \to X$ is the branched cover corresponding to the hyperelliptic involution of $T^2$.  In this case $X$ is the sphere with four marked points.  Rotation of $T^2$ by $\pi$ in one factor is a fiber-preserving homeomorphism homotopic to the identity, but the induced homeomorphism of $X$ acts nontrivially on the marked points and hence is not homotopic to the identity.  Thus this cover fails the Birman--Hilden property.  One can construct a similar example when $S$ is the sphere $S^2$ and $p : S^2 \to X$ is the branched cover induced by a finite order rotation.

%%%
%%%
%%%

\section{Restatement of the main theorem}

We will now give an interpretation of the Birman--Hilden property---hence all three theorems above---in terms of mapping class groups.  Here, the \emph{mapping class group} of a surface is the group of homotopy classes of orientation-preserving homeomorphisms that fix the boundary pointwise and preserve the set of marked points (homotopies must also fix the boundary and preserve the set of marked points).

Let $p: S \to X$ be a covering map of surfaces, possibly branched.   We treat each branch point in $X$ as a marked point, and so homeomorphisms of $X$ are assumed to preserve the set of branch points.  Let $\LMod(X)$ denote the subgroup of the mapping class group $\Mod(X)$ consisting of all elements that have representatives that lift to homeomorphisms of $X$.  This group is called the \emph{liftable mapping class group} of $X$.

Let $\SMod(S)$ denote the subgroup of $\Mod(S)$ consisting of the homotopy classes of all fiber-preserving homeomorphisms.  Here we  emphasize that two homeomorphisms of $S$ are identified in $\SMod(S)$ if they differ by an isotopy that is not necessarily fiber preserving (so that we have a subgroup of $\Mod(S)$).    We also emphasize that preimages of branch points are not marked.  Fiber-preserving homeomorphisms are also called \emph{symmetric homeomorphisms}; these are exactly the lifts of liftable homeomorphisms of $X$.  The group $\SMod(S)$ is called the \emph{symmetric mapping class group} of $S$.

Let $D$ denote the subgroup of $\SMod(S)$ consisting of the homotopy classes of the deck transformations (it is a fact that nontrivial deck transformations represent nontrivial mapping classes).  

\begin{proposition}\label{mcg}
Let $p : S \to X$ be a finite-sheeted branched covering map where $S$ is a hyperbolic surface without boundary.  Then the following are equivalent:
\begin{itemize}
\item $p$ has the Birman--Hilden property,
\item the natural map $\LMod(X) \to \SMod(S)/D$ is injective,
\item the natural map $\SMod(S) \to \LMod(X)$ is well defined, and
\item $\SMod(S)/D \cong \LMod(X)$.
\end{itemize}
\end{proposition}

The proposition is straightforward to prove.  The main content is the equivalence of the first two statements.  The other statements, while useful in practice, are equivalent by rudimentary abstract algebra.  Using the proposition, one obtains several restatements of Theorems~\ref{bh}, \ref{mh}, and~\ref{win} in terms of mapping class groups.

Birman and Hilden also proved that for a regular cover $\SMod(S)$ is the normalizer in $\Mod(S)$ of the deck group $D$ (regarded as a subgroup of $\Mod(S)$), and so we can also write the last statement in Proposition~\ref{mcg} as
\[ 
N_{\Mod(S)}(D) / D \cong \LMod(X).
\]
Birman and Hilden only stated the result about normalizers in the case where the deck group is cyclic.  However, by combining their argument with Kerckhoff's resolution of the Nielsen realization problem \cite{kerckhoff} one obtains the more general version.

There is also a version of Proposition~\ref{mcg} for surfaces with boundary.  Since the mapping class group of a surface with boundary is torsion free, the deck transformations do not represent elements of $\Mod(S_g)$.  And so in this case we can simply replace $D$ with the trivial group.  For example, in the presence of boundary the Birman--Hilden property is equivalent to the statement that $\SMod(S) \cong \LMod(X)$.  This will become especially important in the discussion of braid groups below.

%%%
%%%
%%%

\section{Application to presentations of mapping class groups} 

The original work on the Birman--Hilden theory concerns the case of the hyperelliptic involution and is reported in the 1971 paper \emph{On the mapping class groups of closed surfaces as covering spaces} \cite{BH71}.  We will explain how Theorem~\ref{bh} specializes in this case and helps to give presentations for the associated symmetric mapping class group and the full mapping class group in genus two.

Consider the covering space $S_g \to S_{0,2g+2}$ induced by a hyperelliptic involution of $S_g$.   In general a \emph{hyperelliptic involution} of $S_g$ is a homeomorphism of order two that acts by $-I$ on $H_1(S_g;\Z)$; we remark that the hyperelliptic involution is unique for $S_1$ and $S_2$ but there are infinitely many hyperelliptic involutions of $S_g$ when $g \geq 3$, all conjugate in the group of homeomorphisms.  

Theorem~\ref{bh} and Proposition~\ref{mcg} give an isomorphism
\[ 
\SMod(S_g) / \langle \iota \rangle \cong \LMod(S_{0,2g+2}).
\]
In the special case of the hyperelliptic involution we have $\LMod(S_{0,2g+2}) = \Mod(S_{0,2g+2})$.  Indeed, we can check directly that each half-twist generator for $\Mod(S_{0,2g+2})$ lifts to a Dehn twist in $S_g$.  

In the case $g=2$ we further have 
\[ 
\SMod(S_2) = \Mod(S_2).
\]
In other words, every mapping class of $S_2$ is symmetric with respect to the hyperelliptic involution.  The easiest way to see this is to note that each of the Humphries generators for $\Mod(S_2)$ is a Dehn twist about a curve that is preserved by the hyperelliptic involution.  We thus have the following isomorphism:
\[ \Mod(S_2) / \langle \iota \rangle \cong \Mod(S_{0,6}). \]
Simple presentations for $\Mod(S_{0,n})$ were found by Magnus, and so from his presentation for $\Mod(S_{0,6})$ Birman and Hilden use the above isomorphism to derive the following presentation for $\Mod(S_2)$.  The generators are the Humphries generators for $\Mod(S_2)$, and we denote them by $T_1,\dots, T_5$.  The relations are:  \begin{align*}
[T_i,T_j]&=1 \hspace*{.25in}\text{ for } |i-j|>2 \\
T_iT_{i+1}T_i=T_{i+1}T_i&T_{i+1} \hspace*{.25in}\text{ for } 1 \leq i \leq 4 \\
(T_1T_2T_3T_4T_5)^6&=1 \\
(T_1T_2T_3T_4T_5T_5T_4T_3T_2T_1)^2&=1 \\
\left[T_1T_2T_3T_4T_5T_5T_4T_3T_2T_1,T_1\right]&=1
\end{align*}
The first two relations are the standard braid relations from $B_6$, the next relation describes the kernel of the map $B_6 \to \Mod(S_{0,6})$, and the last two relations come from the two-fold cover: the mapping class
\[ 
T_1T_2T_3T_4T_5T_5T_4T_3T_2T_1
\]
is the hyperelliptic involution.  This presentation is the culmination of a program begun by Bergau and Mennicke \cite{BG}, who approached the problem by studying the surjective homomorphism $B_6 \to \Mod(S_2)$ that factors through the map $\Mod(S_{0,6}) \to \Mod(S_2)$ used here.

Birman used the above presentation to give a normal form for elements of $\Mod(S_2)$ and hence a method for enumerating 3-manifolds of Heegaard genus two \cite{birman}.

As explained by Birman and Hilden, the given presentation for $\Mod(S_2)$ generalizes to a presentation for $\SMod(S_g)$.  The latter presentation has many applications to the study of $\SMod(S_g)$.  It was used by Meyer \cite{meyer72} to show that if a surface bundle over a surface has monodromy in $\SMod(S_g)$ then the signature of the resulting 4-manifold is zero; see also the related work of Endo \cite{endo}.  Endo and Kotschick used the Birman--Hilden presentation to show that the second bounded cohomology of $\SMod(S_g)$ is nontrivial \cite{EK01}.  Also, Kawazumi \cite{kawazumi} used it to understand the low-dimensional cohomology of $\SMod(S_g)$.  

In 1972 Birman and Chillingworth published the paper \emph{On the homeotopy group of a non-orientable surface} \cite{BC}.  There, they determine a generating set for the mapping class group (= homeotopy group) of an arbitrary closed non-orientable surface using similar ideas, namely, they exploit the associated orientation double cover and pass information through the Birman--Hilden theorem from the orientable case.  They also find an explicit finite presentation for the mapping class group of a closed non-orientable surface of genus three, which admits a degree two cover by $S_2$.

One other observation from the 1971 paper is that $\Mod(S_2)$ is both a quotient of and a subgroup of $\Mod(S_{2,6})$.  To realize $\Mod(S_2)$ as a quotient, we consider the map $\Mod(S_{2,6}) \to \Mod(S_2)$ obtained by forgetting the marked points/punctures; this is a special case of the Birman exact sequence studied by Birman in her thesis \cite{birmanthesis}.   And to realize $\Mod(S_2)$ as a subgroup, we use the Birman--Hilden theorem: since every element of $\Mod(S_2)$ has a symmetric representative that preserves the set of preimages of the branch points in $S_{0,6}$ and since isotopies between symmetric homeomorphisms can also be chosen to preserve this set of six points, we obtain the desired inclusion.  Birman and Hilden state that ``the former property is easily understood but the latter much more subtle.'' As mentioned by Mess \cite{mess}, the inclusion $\Mod(S_2) \to \Mod(S_{2,6})$ can be rephrased as describing a multi-section of the universal bundle over moduli space in genus two.

%%%
%%%
%%%

\section{More applications to the genus two mapping class group} 

In the previous section we saw how the Birman--Hilden theory allows us to transport knowledge about the mapping class group of a punctured sphere to the mapping class group of a surface of genus two.  As the former are closely related to braid groups, we can often push results about braid groups to the mapping class group.  Almost every result about mapping class groups that is special to genus two is proved in this way.  

A prime example of this is the result of Bigelow--Budney \cite{BB} and Korkmaz \cite{KL} which states that $\Mod(S_2)$ is linear, that is, $\Mod(S_2)$ admits a faithful representation into $\GL_N(\C)$ for some $N$.  Bigelow and Krammer independently proved that braid groups were linear, and so the main work is to derive from this the linearity of $\Mod(S_{0,n})$.  They then use the isomorphism $\Mod(S_2) / \langle \iota \rangle \cong \Mod(S_{0,6})$ to push the linearity up to $\Mod(S_2)$.

A second example is from the thesis of Whittlesey, published in 2000.  She showed that $\Mod(S_2)$ contains a normal subgroup where every nontrivial element is pseudo-Anosov \cite{whittlesey}.  The starting point is to consider the Brunnian subgroup of $\Mod(S_{0,6})$.  This is  the intersection of the kernels of the six forgetful maps $\Mod(S_{0,6}) \to \Mod(S_{0,5})$, so it is obviously normal in $\Mod(S_{0,6})$.  She shows that all nontrivial elements of this group are pseudo-Anosov and proves that the preimage in $\Mod(S_2)$ has a finite-index subgroup with the desired properties.

We give one more example.  In the 1980s, before the work of Bigelow and Krammer, Vaughan Jones discovered a representation of the braid group defined in terms of Hecke algebras \cite{jones}.  As in the work of  Bigelow--Budney and Korkmaz, one can then derive a representation of $\Mod(S_{0,2g+2})$ and then---using the Birman--Hilden theory---of $\SMod(S_g)$.  When $g=2$ we thus obtain a representation of $\Mod(S_2)$ to $\GL_5(\Z[t,t^{-1}])$.  This representation was used by Humphries \cite{humphries} to show that the normal closure in $\Mod(S_2)$ of the $k$th power of a Dehn twist about a nonseparating curve has finite index if and only if $|k| \leq 3$.

There are many other examples, such as the computation of the asymptotic dimension of $\Mod(S_2)$ by Bell and Fujiwara \cite{BF} and the determination of the minimal dilatation in $\Mod(S_2)$ by Cho and Ham \cite{CH}; the list goes on, but so must we.

%%%
%%%
%%%

\section{Application to representations of the braid group}\label{rep}

The Birman--Hilden theorem gives an important relationship between braid groups and mapping class groups.  This is probably the most oft-used application of their results.  

Let $S_g^1$ the orientable surface of genus $g$ with one boundary component and let $D_{2g+1}$ denote the closed disk with $2g+1$ marked points in the interior.  Consider the covering space $S_g^1 \to D_{2g+1}$ induced by a hyperelliptic involution of $S_g^1$.  It is well known that $\Mod(D_{2g+1})$ is isomorphic to the braid group $B_{2g+1}$.  As in the closed case, it is not hard to see that $\LMod(D_{2g+1}) = \Mod(D_{2g+1})$ (again, each of the standard generators for $B_{2g+1}$ lifts to a Dehn twist).

One is thus tempted to conclude that $\SMod(S_g^1) / \langle \iota \rangle \cong B_{2g+1}$.  But this is not the right statement, since $\iota$ does not represent an element of $\Mod(S_g^1)$.  Indeed, for surfaces with boundary we insist that homeomorphisms and homotopies fix the boundary pointwise (otherwise we would not have the isomorphism $\Mod(D_{2g+1}) \cong B_{2g+1}$).  Therefore, the correct isomorphism is:
\[ \SMod(S_g^1) \cong B_{2g+1} \]
The most salient aspect of this isomorphism is that there is an injective homomorphism
\[ B_{2g+1} \to \Mod(S_g^1). \]
The injectivity here is sometimes attributed to Perron--Vannier \cite{PV}.  It is possible that they were the first to observe this consequence of the Birman--Hilden theorem but the only nontrivial step is the Birman--Hilden theorem.  

In the case of $g=1$ the representation of $B_3$ is onto $\Mod(S_1^1)$, and so 
\[ \Mod(S_1^1) \cong B_3. \]
Similarly we have
\[ \Mod(S_1^2) \cong B_4 \times \Z. \]
The point here is that $B_4$ surjects onto $\SMod(S_1^2)$ and the latter is almost isomorphic to $\Mod(S_1^2)$; the extra $\Z$ comes from the Dehn twist about a single boundary component.

One reason that the embedding of $B_{2g+1}$ in $\Mod(S_g^1)$ is so important is that if we compose with the standard symplectic representation
\[ \Mod(S_g^1) \to \Sp_{2g}(\Z) \]
then we obtain a representation of the braid group
\[ B_{2g+1} \to \Sp_{2g}(\Z) .\]
This representation is called the standard symplectic representation of the braid group.  It is also called the \emph{integral Burau representation} because it is the only integral specialization of the Burau representation besides the permutation representation.  The symplectic representation is obtained by specializing the Burau representation at $t=-1$, while the permutation representation is obtained by taking $t=1$.

The image of the integral Burau representation has finite index in the symplectic group: it is an extension of the level two symplectic group by the symmetric group on $2g+1$ letters.  The projection onto the symmetric group factor is the standard symmetric group representation of the braid group.  See A'Campo's paper \cite{acampo} for details.

The kernel of the integral Burau representation is known as the hyperelliptic Torelli group.  This group is well studied, as it describes the fundamental group of the branch locus of the period mapping from Teichm\"uller space to the Siegel upper half-space; see, for instance, the paper by Brendle, Putman, and the first author of this article \cite{BMP} and the references therein.

There are plenty of variations on the given representation.  Most important is that if we take a surface with two boundary components $S_g^2$ and choose a hyperelliptic involution, that is, an order two homeomorphism that acts by $-I$ on the first homology of the surface, then the quotient is $D_{2g+2}$ and so we obtain an isomorphism:
\[ \SMod(S_g^2) \cong B_{2g+2}.\]
Also, since the inclusions $S_g^1 \to S_{g+1}$ and $S_g^2 \to S_{g+1}$ induce injections $\SMod(S_g^1) \to \Mod(S_{g+1})$ and $\SMod(S_g^2) \to \Mod(S_{g+1})$ we obtain embeddings of braid groups into mapping class groups of closed surfaces.

In the 1971 paper Birman and Hilden discuss the connection with representations of the braid group.  They point out the related fact that $B_{2g+2}$ surjects onto $\SMod(S_g)$ (this follows immediately from their presentation for the latter).  In the special case $g=1$ this becomes the classical fact that $B_4$ surjects onto $\SMod(S_1) = \Mod(S_1) \cong \SL_2(\Z)$.  We can also derive this fact from our isomorphism $\Mod(S_1^1) \cong B_3$, the famous surjection $B_4 \to B_3$, and the surjection $\Mod(S_1^1) \to \Mod(S_1)$ obtained by capping the boundary.

One useful application of the embeddings of braid groups in mapping class groups is that we can often transport relations from the former to the latter.  In fact, almost all of the widely-used relations in the mapping class group have interpretations in terms of braids.  This is especially true in the theory of Lefschetz fibrations; see for instance the work of Korkmaz \cite{korkmaz01} and Hamada \cite{hamada} and of Baykur and Van Horn-Morris \cite{BV}.

%%%
%%%
%%%

\section{Application to a question of Magnus}

The last application we will explain is beautiful and unexpected.  It is the resolution of a seemingly unrelated question of Magnus about braid groups.  

As mentioned in the previous section, the braid group $B_n$ is isomorphic to the mapping class group of a disk $D_n$ with $n$ marked points.  Let us write $D_n^\circ$ for the surface obtained by removing from $D_n$ the marked points.  There is then a natural action of $B_n$ on $\pi_1(D_n^\circ)$ (with base point on the boundary).  The latter is isomorphic to the free group $F_n$ on $n$ letters.  Basic algebraic topology tells us that this action is faithful.  In other words, we have an injective homomorphism:
\[ 
B_n \to \Aut(F_n).
\]
This is a fruitful way to view the braid group; for instance, since the word problem in $\Aut(F_n)$ is easily solvable, this gives a solution to the word problem for $B_n$.  

Let $F_{n,k}$ denote the normal closure in $F_n$ of the elements $x_1^k,\dots,x_n^k$.  The quotient $F_n/F_{n,k}$ is isomorphic to the $n$-fold free product $\Z/k\Z \ast \cdots \ast \Z/k\Z$.  Since the elements of $B_n$ preserve the set of conjugacy classes $\{[x_1],\dots,[x_n]\}$, there is an induced homomorphism
\[
B_n \to \Aut(F_n/F_{n,k}).
\]
Let $B_{n,k}$ denote the image of $B_n$ under this map.  Magnus asked:
\begin{quote}
\emph{Is $B_n$ isomorphic to $B_{n,k}$?}
\end{quote}
In other words, is the map $B_n \to \Aut(F_n/F_{n,k})$ injective?

In their \emph{Annals paper}, Birman and Hilden answer Magnus' question in the affirmative.  Here is the idea.  Let $H_{n,k}$ denote the kernel of the map
\[ F_n \to \Z/k\Z \]
where each generator of $F_n$ maps to 1.  The covering space of $D_n^\circ$ corresponding to $H_{n,k}$ is a $k$-fold cyclic cover $S^\circ$.  If we consider a small neighborhood of one of the punctures in $D_n^\circ$, the induced covering map is equivalent to the connected $k$-fold covering space of $\C \setminus \{0\}$ over itself (i.e. the one induced by $z \mapsto z^k$).  As such, we can ``plug in'' to $S^\circ$ a total of $n$ points in order to obtain a surface $S$ and a cyclic branched cover $S \to D_n$.   The fundamental group of $S^\circ$ is $H_{n,k}$ by definition.  It follows from Van Kampen's theorem that $\pi_1(S) \cong H_{n,k}/F_{n,k}$.  Indeed, a simple loop around a puncture in $S^\circ$ projects to a loop in $D_n^\circ$ that circles the corresponding puncture $k$ times.

As in the case of the hyperelliptic involution, we can check directly that each element of $B_n$ lifts to a fiber-preserving homeomorphism of $S$.  Therefore, to answer Magnus' question in the affirmative it is enough to check that the map $B_n \to \Aut \pi_1(S)$ is injective.  Suppose $b \in B_n$ lies in the kernel.  Then the corresponding fiber-preserving homeomorphism of $S$ is homotopic, hence isotopic, to the identity.  By the Birman--Hilden theorem (the version for surfaces with boundary), $b$ is trivial, and we are done.

Bacardit and Dicks \cite{BD} give a purely algebraic treatment of Magnus' question; they credit the argument to Crisp and Paris \cite{CP}.  Another algebraic argument for the case of even $k$ was given by D.L. Johnson \cite{johnson}.   Yet another combinatorial proof was given in 1992 by Kr\"uger \cite{kruger}.

%%%
%%%
%%%

\section{A famous (but false) proof of the Birman--Hilden theorem}

When confronted with the Birman--Hilden theorem, one might be tempted to quickly offer the following easy proof: given the branched cover $p : S \to X$, the fiber-preserving homeomorphism $f : S \to S$, the corresponding homeomorphism $\bar f : X \to X$, and an isotopy $H : S \times I \to S$ from $f$ to the identity map, we can consider the composition $p \circ H$, which gives a homotopy from $\bar f$ to the identity.  Then, since homotopic homeomorphisms of a surface are isotopic, there is an isotopy from $\bar f$ to the identity, and this isotopy lifts to a fiber-preserving isotopy from $f$ to the identity.  Quod erat demonstrandum.

This probably sounds convincing, but there are two problems.  First of all, the composition $p \circ H$ is really a homotopy between $p \circ f$ and $p$ which are maps from $S$ to $X$; since $H$ is not fiber preserving, there is no way to convert this to a well-defined homotopy between maps $X \to X$.  The second problem is that $H$ might send points that are not preimages of branch points to preimages of branch points; so even if we could project the isotopy, we would not obtain a homotopy of $X$ that respects the marked points.

In the next two sections will will outline proofs of the Birman--Hilden theorem in various cases.  The reader should keep in mind the subtleties uncovered by this false proof.

\section{The unbranched (= easy) case}

Before getting to the proof of the Birman--Hilden theorem, we will warm up with the case of unbranched covers.  This case is much simpler, as all of the subtlety of the Birman--Hilden theorem lies in the branch points.  Still the proof is nontrivial, and later we will prove the more general case by reducing to the unbranched case.

In 1972 Birman and Hilden published the paper \emph{Lifting and projecting homeomorphisms} \cite{BH72}, which gives a quick proof of Theorem~\ref{bh} in the case of regular unbranched covers.  Following along the same lines, Aramayona, Leininger, and Souto generalized their proof to the case of arbitrary (possibly irregular) unbranched covers \cite{ALS}.  We will now explain their proof.  

Let $p : S \to X$ be an unbranched covering space of surfaces, and let $f : S \to S$ be a fiber-preserving homeomorphism that is isotopic to the identity.  Without loss of generality, we may assume that $f$ has a fixed point.  Indeed, if $f$ does not fix some point $x$, then we can push $p(x)$ in $X$ by an ambient isotopy, and lift this isotopy to $S$ until $x$ is fixed.  As a consequence, $f$ induces a well-defined action $f_\star$ on $\pi_1(S)$.  Since $f$ is isotopic to the identity, $f_\star$ is the identity.  If $\bar f$ is the corresponding homeomorphism of $X$, then it follows that $\bar f_\star$ is the identity on the finite-index subgroup $p_\star(\pi_1(S))$ of $\pi_1(X)$.  From this, plus the fact that roots are unique in $\pi_1(X)$, we conclude that $\bar f_\star$ is the identity.  By basic algebraic topology, $\bar f$ is homotopic to the identity, and hence it is isotopic to the identity, which implies that $f$ is isotopic to the identity through fiber-preserving homeomorphisms, as desired.

%%%
%%%
%%%

\section{Three (correct) proofs of the Birman--Hilden theorem}

In this section we present sketches of the proofs of all three versions of the Birman--Hilden theorem given in Section~\ref{state}.  We begin with the original proof by Birman and Hilden, which is a direct attack using algebraic and geometric topology.  Then we explain the proof of MacLachlan and Harvey's Teichm\"uller theoretic approach, and finally the combinatorial topology approach of the second author, which gives a further generalization.  

\p{The Birman--Hilden proof: Algebraic and geometric topology} As in the statement of Theorem~\ref{bh}, let $p : S \to X$ be a regular branched covering space  where $S$ is a hyperbolic surface.  As in Theorem~1 of the \emph{Annals} paper by Birman and Hilden, we make the additional assumption here that each deck transformation for this cover fixes each preimage of each branch point in $X$.  Theorem~\ref{bh} will follow easily from this special case.  Let $f$ be a fiber-preserving homeomorphism of $S$ and assume that $f$ is isotopic to the identity.  

Let $x$ be the preimage in $S$ of some branch point in $X$.  The first key claim is that $f(x)=x$ \cite[Lemma 1.3]{BH73}.  Thus if we take the isotopy $H$ from $f$ to the identity and restrict it to $x$, we obtain an element $\alpha$ of $\pi_1(S,x)$.  Birman and Hilden argue that $\alpha$ must be the trivial element.  The idea is to argue that $\alpha$ is fixed by each deck transformation (this makes sense since the deck transformations fix $x$), and then to argue that the only element of $\pi_1(S)$ fixed by a nontrivial deck transformation is the trivial one (to see this, regard $\alpha$ as an isometry of the universal cover $\mathbb{H}^2$ and regard a deck transformation as a rotation of $\mathbb{H}^2$).

Since $\alpha$ is trivial, we can deform it to the trivial loop, and by extension we can deform the isotopy $H$ to another isotopy that fixes $x$ throughout.  Proceeding inductively, Birman and Hilden argue that $H$ can be deformed so that it fixes all preimages of branch points throughout the isotopy.  At this point, by deleting branch points in $X$ and their preimages in $S$, we reduce to the unbranched case.

Finally, to prove their Theorem~2, which treats the case of solvable covers, Birman and Hilden reduce it to Theorem~1 by factoring any solvable cover into a sequence of cyclic covers of prime order. Such a cover must satisfy the hypotheses of Theorem~1.

It would be interesting to use the Birman--Hilden approach to prove the more general theorem of Winarski.  There is a paper by  Zieschang from 1973 that uses similar reasoning to Birman and Hilden and recovers the result of MacLachlan and Harvey \cite{zieschang73}.  

\p{MacLachlan and Harvey's proof: Teichm\"uller theory} We now explain the approach of MacLachlan and Harvey.  Let $p : S \to X$ be a regular branched covering space where $S$ is a hyperbolic surface.  We will give MacLachlan and Harvey's argument for Theorem~\ref{mh} and at the same time explain why the argument gives the more general result of Theorem~\ref{win}.

The mapping class group $\Mod(S)$ acts on the Teichm\"uller space $\Teich(S)$, the space of isotopy classes of complex structures on $S$ (or conformal structures on $S$, or hyperbolic structures on $S$, or algebraic structures on $S$).  Let $X^\circ$ denote the complement in $X$ of the set of branch points.  There is a map $\Xi : \Teich(X^\circ) \to \Teich(S)$ defined by lifting complex structures through the covering map $p$ (one must apply the removable singularity theorem to extend over the preimages of the branch points).  

The key point in the proof is that $\Xi$ is injective.  One way to see this is to observe that Teichm\"uller geodesics in $\Teich(X^\circ)$ map to Teichm\"uller geodesics in $\Teich(S)$ of the same length.  Indeed, the only way this could fail would be if we had a Teichm\"uller geodesic in $\Teich(X^\circ)$ where the corresponding quadratic differential had a simple pole (= 1-pronged singularity) at a branch point and some pre-image of that branch point was unramified (1-pronged singularities are only allowed at marked points, and preimages of branch points are not marked).  This is why the most natural setting for this argument is that of Theorem~\ref{win}, namely, where $p$ is fully ramified.

Let $Y$ denote the image of $\Xi$.  The symmetric mapping class group $\SMod(S)$ acts on $Y$ and the kernel of this action is nothing other than $D$.  The liftable mapping class group $\LMod(X)$ acts faithfully on $\Teich(X^\circ)$ and hence---as $\Xi$ is injective---it also acts faithfully on $Y$.  It follows immediately from the definitions that the images of $\SMod(S)$ and $\LMod(X)$ in the group of automorphisms of $Y$ are equal.  It follows that $\SMod(S)/D$ is isomorphic to $\LMod(X)$, as desired.

\p{Winarski's proof: Combinatorial topology}   Let $p : S \to X$ be a fully ramified branched covering space where $S$ is a hyperbolic surface.  To prove Theorem~\ref{win} we will show that $\Phi : \LMod(X) \to \SMod(S)/D$ is injective.  

Suppose $f \in \LMod(X)$ lies in the kernel of $\Phi$.  Let $\varphi$ be a representative of $f$.  Since $\Phi(f)$ is trivial we can choose a lift $\tilde \varphi : S \to S$  that is isotopic to the identity; thus $\tilde \varphi$ fixes the isotopy class of every simple closed curve in $S$.  The main claim is that $\varphi$ fixes the isotopy class of every simple closed curve in $X$.  From this, it follows that $f$ has finite order in $\LMod(X)$.  Since $\ker(\Phi)$ is torsion free \cite[Prop 4.2]{winarski}, the theorem will follow.

So let us set about the claim.  Let $c$ be a simple closed curve in $X$, and let $\tilde c$ be its preimage in $S$.  By assumption $\tilde \varphi(\tilde c)$ is isotopic to $\tilde c$ and we would like to leverage this to show $\varphi(c)$ is isotopic to $c$.  There are two stages to the argument: first dealing with the case where $\varphi(c)$ and $c$ are disjoint, and then in the case where they are not disjoint we reduce to the disjoint case.

If $\varphi(c)$ and $c$ are disjoint, then $\tilde \varphi( \tilde c)$ and $\tilde c$ are disjoint.  Since the latter are isotopic, they co-bound a collection of annuli $A_1,\dots,A_n$.  Then, since orbifold Euler characteristic is multiplicative under covers, we can conclude that $p(\cup A_i)$ is an annulus with no branch points (branch points decrease the orbifold Euler characteristic), and so $c$ and $\varphi(c)$ are isotopic.

We now deal with the second stage, where $\varphi(c)$ and $c$ are not disjoint.  In this case $\tilde \varphi(\tilde c)$ and $\tilde c$ are not disjoint either, but by our  assumptions they are isotopic in $S$.  Therefore, $\tilde \varphi(\tilde c)$ and $\tilde c$ bound at least one bigon.

Consider an innermost such bigon $B$.  Since $B$ is innermost, $p(B)$ is an innermost bigon bounded by $\varphi(c)$ and $c$ in $X$ (the fact that $B$ is innermost implies that $p|B$ is injective).  If there were a branch point in $p(B)$ then since $p$ is fully ramified, this would imply that $B$ was a $2k$-gon with $k > 1$, a contradiction.  Thus, we can apply an isotopy to remove the bigon $p(B)$ and by induction we reduce to the case where $\varphi(c)$ and $c$ are disjoint.

For an exposition of Winarski's proof in the case of the hyperelliptic involution, see the book by Farb and the first author of this article \cite{primer}.

%%%
%%%
%%%

\begin{center}
\begin{figure}
\includegraphics{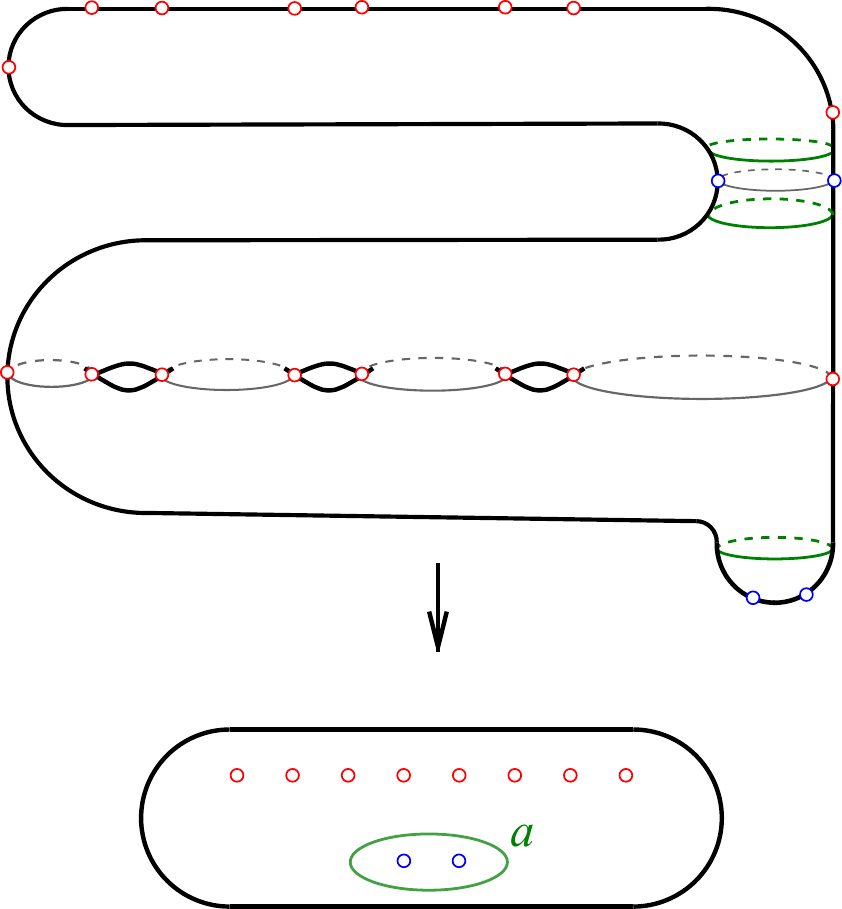}
\caption{The simple threefold cover of $S_g$ over the sphere.  The gray curves divide $S_g$ into three regions, each serving as a ``fundamental domain'' for the cover.  The preimages in $S_g$ of the branch points in $S^2$ are shown as dots but are treated as unmarked points in $S_g$.}
\label{fuller}
\end{figure}
\end{center}

\section{Open questions and new directions}

One of the most striking aspects of the Birman--Hilden story is the breadth of open problems related to the theory and the constant discovery of related directions.  We  have already mentioned a number of questions arising from the theory.  
Perhaps the most obvious open problem is the following.

\begin{question}
Which branched covers of surfaces have the Birman--Hilden property?
\end{question}

Based on the discussion in Section~\ref{state} above, one might hope that all branched coverings---at least where the cover is a hyperbolic surface---have the Birman--Hilden property.  However, this is not true.  Consider, for instance, the simple threefold cover $p : S_g \to X$, where $X$ is the sphere with $2g+4$ branch points (this cover is unique).  As shown in Figure~\ref{fuller} we can find an essential curve $a$ in $X$ whose preimage  in $S_g$ is a union of three homotopically trivial simple closed curves.   It follows that the Dehn twist $T_a$ lies in the kernel of the map $\LMod(X) \to \SMod(S_g)$ and so $p$ does not have the Birman--Hilden property.  See Fuller's paper for further discussion of this example and the relationship to Lefschetz fibrations \cite{fuller}.  

Berstein and Edmonds generalized this example by showing that no simple cover of degree at least three over the sphere has the Birman--Hilden property \cite{BE}, and Winarski further generalized this by proving that no simple cover of degree at least three over any surface has the Birman--Hilden property \cite{winarski}.

Having accepted the fact that not all covers have the Birman--Hilden property, one's second hope might be that a cover has the Birman--Hilden property if and only if it is fully ramified.  However, this is also false.  Chris Leininger \cite{leininger} has explained to us how to construct a counterexample using the following steps.  First, let $S$ be a surface and let $z$ be a marked point in $S$.  Let $p :  \tilde S \to S$ be a characteristic cover of $S$ and let $\tilde z$ be one point of the full preimage $p^{-1}(z)$.  Ivanov and McCarthy \cite{IM} observed that there is an injective homomorphism $\Mod(S,z) \to \Mod(\tilde S,p^{-1}(z))$ where for each element of $\Mod(S,z)$, we choose the lift to $\tilde S$ that fixes $\tilde z$.  Aramayona--Leininger--Souto \cite{ALS} proved that the composition of the Ivanov--McCarthy homomorphism with the forgetful map $\Mod(\tilde S,p^{-1}(z)) \to \Mod(\tilde S, \tilde z)$ is injective.  If we then take a regular branched cover $S' \to \tilde S$ with branch locus $\tilde z$, the resulting cover $S' \to S$ is not fully ramified but it has the Birman--Hilden property.

Here is another basic question.

\begin{question}\label{smodq}
For which cyclic branched covers of $S_g$ over the sphere is $\SMod(S_g)$ equal to a proper subgroup of $\Mod(S_g)$?  When is it finite index?
\end{question}

Theorem 5 in the \emph{Annals} paper by Birman and Hilden states for a cyclic branched cover $S \to X$ over the sphere we have $\LMod(X) =\Mod(X)$.  Counterexamples to this theorem were recently discovered by Ghaswala and the second author (see the erratum \cite{erratum}), who wrote a paper \cite{GW} classifying exactly which branched covers over the sphere have $\LMod(X)=\Mod(X)$.  Theorem 6 in the paper by Birman and Hilden states that for a cyclic branched cover of $S_g$ over the sphere with $g \geq 3$ the group $\SMod(S_g)$ is a proper subgroup of $\Mod(S_g)$.  The proof uses their Theorem 5, so Question~\ref{smodq} should be considered an open question.   Of course this question can be generalized to other base surfaces besides the sphere and other types of covers.  For simple branched covers over the sphere (which, as above, do not have the Birman--Hilden property) Berstein and Edmonds \cite{BE} proved that $\SMod(S_g)$ is equal to $\Mod(S_g)$.  

\medskip

We can also ask about the Birman--Hilden theory for orbifolds and 3-manifolds.

\begin{question}
Which covering spaces of two-dimensional orbifolds have the Birman--Hilden property?
\end{question}

Earle proved some Birman--Hilden-type results for orbifolds in his recent paper \cite{earle}, which he describes as a sequel to his 1971 paper \emph{On the moduli of closed Riemann surfaces with symmetries} \cite{earle71}.

\begin{question}
Which covering spaces of 3-manifolds enjoy the Birman--Hilden property?
\end{question}

Vogt proved that certain regular unbranched covers of certain Seifert-fibered 3-manifolds have the Birman--Hilden property \cite{vogt}.  He also explains the connection to understanding foliations in codimension two, specifically for foliations of closed 5-manifolds by Seifert 3-manifolds.

A specific 3-manifold worth investigating is the connect sum of $n$ copies of $S^2 \times S^1$; call it $M_n$.  The outer automorphism group of the free group $F_n$ is a finite quotient of the mapping class group of $M_n$.  Therefore, one might obtain a version of the Birman--Hilden theory for the outer automorphism group of $F_n$ by developing a Birman--Hilden theory for $M_n$.  

\begin{question}
Does $M_n$ enjoy the Birman--Hilden property?  If so, does this give a Birman--Hilden theory for free groups?
\end{question}

For example, consider the hyperelliptic involution $\sigma$ of $F_n$, the outer automorphism that (has a representative that) inverts each generator of $F_n$.  This automorphism is realized by the homeomorphism of $M_n$ that reverses each $S^1$-factor.  The resulting quotient of $M_n$ is the 3-sphere with branch locus the $(n+1)$-component unlink.  This is in consonance with the fact that the centralizer of $\sigma$ in the outer automorphism group of $F_n$ is the palindromic subgroup and that the latter is closely related to the configuration space of unlinks in $S^3$; see the paper by Collins \cite{collins}.

\medskip

Next, there are many questions about the hyperelliptic Torelli group and its generalizations.  As discussed in Section~\ref{rep}, the hyperelliptic Torelli group is the kernel of the integral Burau representation of the braid group.  With Brendle and Putman, the first author of this article proved \cite{BMP}  that this group is generated by the squares of Dehn twists about curves that surround an odd number of marked points in the disk $D_n$.

\begin{question}
Is the hyperelliptic Torelli group finitely generated?  Is it finitely presented?  Does it have finitely generated abelianization?
\end{question}

There are many variants of this question.  By changing the branched cover $S \to D_n$, we obtain many other representations of (the liftable  subgroups of) the braid group.  Each representation gives rise to its own Torelli group.  Except for the hyperelliptic involution case, very little is known.  One set of covers to consider are the superelliptic covers studied by Ghaswala and the second author of this article \cite{GW17}.

Another aspect of this question is to determine the images of the braid groups in $\Sp_{2N}(\Z)$ under the various representations of (finite index subgroups of) the braid group arising from various covers $S \to D_n$.  By work of McMullen \cite{mcmullen} and Venkatarmana \cite{venkataramana}, it is known that when the degree of the cover is at least three and $n$ is more than twice the degree, the image has finite index in the centralizer of the image of the deck group.

\begin{question}
For which covers $S \to D_n$ does the associated representation of the braid group have finite index in the centralizer of the image of the deck group?
\end{question}

There are still many aspects to the Birman--Hilden theory that we have not touched upon.  Ellenberg and McReynolds \cite{EM} used the theory to prove that every algebraic curve over $  \mathbb{\bar Q}$ is birationally equivalent over $\C$ to a Teichm\"uller curve.  Nikolaev  \cite{nikolaev} uses the embedding of the braid group into the mapping class group to give cluster algebraic representations of braid groups.  Kordek applies the aforementioned result of Ghaswala and the second author of this article to deduce information about the Picard groups of various moduli spaces of Riemann surfaces \cite{kordek}.  A Google search for ``Birman--Hilden'' yields a seemingly endless supply of applications and connections (the \emph{Annals} paper has 139 citations on Google Scholar at the time of this writing).  We hope that the reader is inspired to learn more about these connections and pursue their own developments of the theory.

\p{Acknowledgments.} The authors would like to thank John Etnyre, Tyrone Ghaswala, Allen Hatcher, and Chris Leininger for helpful conversations.

\bibliographystyle{plain}
\bibliography{bh}

\end{document}